\documentclass[11pt]{article}
\usepackage{mathrsfs}
\usepackage{amsfonts}
\usepackage{amssymb}
\usepackage{amsmath}
\usepackage[all]{xy}

\def\Ker{\mathop{\rm Ker}\nolimits}
\def\Hom{\mathop{\rm Hom}\nolimits}
\def\Coker{\mathop{\rm Coker}\nolimits}
\def\Mod{\mathop{\rm Mod}\nolimits}
\def\mod{\mathop{\rm mod}\nolimits}
\def\pd{\mathop{\rm pd}\nolimits}
\def\Gpd{\mathop{\rm Gpd}\nolimits}

\def\Im{\mathop{\rm Im}\nolimits}
\def\Tr{\mathop{\rm Tr}\nolimits}
\def\Ext{\mathop{\rm Ext}\nolimits}

\def\Mod{\mathop{\rm Mod}\nolimits}
\def\mod{\mathop{\rm mod}\nolimits}
\title{\Large \bf Gorenstein Syzygy Modules
}
\author{Chonghui Huang$^{1,2}$, \ Zhaoyong Huang$^{2}$\thanks{{\it E-mail address}: huangzy@nju.edu.cn}
\\{\footnotesize $^{1}${\it Research Institute of Mathematics, University of South China, Hengyang
421001, P.R. China;}} \\ {\footnotesize $^{2}${\it
Department of Mathematics, Nanjing University,
Nanjing 210093, P.R. China}}}
\date{}
\begin{document}
\baselineskip=18pt \maketitle

\begin{abstract}
For any ring $R$ and any positive integer $n$, we prove that a left
$R$-module is a Gorenstein $n$-syzygy if and only if it is an
$n$-syzygy. Over a left and right Noetherian ring, we introduce the
notion of the Gorenstein transpose of finitely generated modules. We
prove that a module $M\in \mod R^{op}$ is a Gorenstein transpose of
a module $A\in \mod R$ if and only if $M$ can be embedded into a
transpose of $A$ with the cokernel Gorenstein projective. Some
applications of this result are given.
\end{abstract}

\vspace{0.5cm}

\centerline{\bf 1. Introduction }

\vspace{0.2cm}

Throughout this paper, $R$ is an associative ring with identity and
$\Mod R$ is the category of left $R$-modules.

In classical homological algebra, the notion of finitely generated
projective modules is an important and fundamental research object.
As a generalization of this notion, Auslander and Bridger introduced
in [AB] the notion of finitely generated modules of Gorenstein
dimension zero over a left and right Noetherian ring. Over a general
ring, Enochs and Jenda introduced in [EJ1] the notion of Gorenstein
projective modules (not necessarily finitely generated). It is well
known that these two notions coincide for finitely generated modules
over a left and right Noetherian ring. In particular, Gorenstein
projective modules share many nice properties of projective modules
(e.g. [AB, C, CFH, CI, EJ1, EJ2, H]).

The notion of a syzygy module was defined via the projective
resolution of modules as follows. For a positive integer $n$, a
module $A\in \Mod R$ is called an {\it $n$-syzygy module} (of $M$)
if there exists an exact sequence $0\to A \to P_{n-1} \to \cdots \to
P_1 \to P_0\to M \to 0$ in $\Mod R$ with all $P_i$ projective.
Analogously, we call $A$ a {\it Gorenstein $n$-syzygy module} (of
$M$) if there exists an exact sequence $0\to A \to G_{n-1} \to
\cdots \to G_1 \to G_0\to M \to 0$ in $\Mod R$ with all $G_i$
Gorenstein projective. It is trivial that an $n$-syzygy module is
Gorenstein $n$-syzygy. In Section 2,
our main result is that for every $n \geq 1$, a Gorenstein
$n$-syzygy module is $n$-syzygy. The following auxiliary
proposition plays a crucial role in proving this main result. Let
$0\to A\to G_1\buildrel {f} \over\to G_0\to M\to 0$ be an exact
sequence in $\Mod R$ with $G_0$ and $G_1$ Gorenstein projective.
Then we have the following exact sequences $0\to A\to P \to G\to
M\to 0$ and $0\to A\to H \to Q\to M\to 0$ in $\Mod R$ with $P$, $Q$
projective and $G$, $H$ Gorenstein projective.

In Section 3, for a left and right Noetherian ring $R$ and a
finitely generated left $R$-module $A$, we introduce the notion of
the Gorenstein transpose of $A$, which is a Gorenstein version of
that of the transpose of $A$. We establish a relation between a
Gorenstein transpose of a module and a transpose of the same module.
We prove that a finitely generated right $R$-module $M$ is a
Gorenstein transpose of a finitely generated left $R$-module $A$ if
and only if $M$ can be embedded into a transpose of $A$ with the
cokernel Gorenstein projective. Then we give some applications of
this result: (1) The direct sum of a
finitely generated Gorenstein projective right $R$-module and a
transpose of a finitely generated left $R$-module $A$ is a
Gorenstein transpose of $A$. (2) For any Gorenstein transpose
and any transpose of a finitely generated left $R$-module, one of them is
$n$-torsionfree if and only if so is the other. (3) A finitely
generated left $R$-module with Gorenstein projective dimension $n$
is a double Gorenstein transpose of a finitely generated left
$R$-module with projective dimension $n$.

\vspace{0.5cm}

\centerline{\bf 2. Gorenstein syzygy modules}

\vspace{0.2cm}

Recall from [EJ1] a module $G \in \Mod R$ is called {\it Gorenstein
projective} if there exists an exact sequence in $\Mod R$:
$$\cdots \to P_1 \to P_0 \to P^0 \to P^1\to \cdots ,$$
such that: (1) All $P_i$ and $P^i$ are projective; (2) After
applying the functor $\Hom _{R}(\ , P)$ the sequence is still exact
for any projective module $P \in \Mod R$; and (3) $G\cong \Im(P_0\to
P^0)$. Let $M$ be a module in $\Mod R$. The {\it Gorenstein
projective dimension} of $M$, denoted by $\Gpd _R(M)$, is defined as
inf$\{ n|$for any exact sequence $0 \to G_{n} \to \cdots \to
G_{1} \to G_{0} \to M \to 0$ in $\Mod R$ with all $G_i$ Gorenstein
projective$\}$. We have $\Gpd _R(M)\geq 0$ and we set $\Gpd _R(M)$ infinity if no such integer
exists (see [EJ1 or H]).

\vspace{0.2cm}

{\bf Lemma 2.1.} {\it Let $0\to M_3\to M_2\to M_1\to 0$ be an exact
sequence in $\Mod R$ with $M_3\neq 0$. If $M_1$ is Gorenstein
projective, then $\Gpd _{R}(M_3)=\Gpd _{R}(M_2)$.}

\vspace{0.2cm}

{\it Proof.} By [H, Theorems 2.24 and 2.20], it is easy to get the
assertion. \hfill $\square$

\vspace{0.2cm}

The following result plays a crucial role in this paper.

\vspace{0.2cm}

{\bf Proposition 2.2.} {\it Let $0\to A\to G_1\buildrel {f} \over\to
G_0\to M\to 0$ be an exact sequence in $\Mod R$ with $G_0$ and $G_1$
Gorenstein projective. Then we have the following exact sequences:
$$0\to A\to P \to G\to M\to 0,$$ and
$$0\to A\to H \to Q\to M\to 0,$$
 in
$\Mod R$ with $P$, $Q$ projective and $G$, $H$ Gorenstein
projective.}

\vspace{0.2cm}

{\it Proof}. Because $G_1$ is Gorenstein projective, there exists an
exact sequence $0\to G_1 \to P \to G_2\to 0$ in $\Mod R$ with $P$
projective and $G_2$ Gorenstein projective. Then we have the
following push-out diagram:
$$\xymatrix{
& & 0 \ar[d] &0 \ar[d] &   \\
0 \ar[r]  & A \ar@{=}[d]  \ar[r] & G_1 \ar[d]
\ar[r]  &\Im f \ar[d] \ar[r] & 0 \\
0 \ar[r] & A \ar[r] & P \ar[d] \ar[r] & B \ar[d] \ar[r] & 0  \\
&  & G_2 \ar[d]
\ar@{=}[r] & G_2\ar[d]  &   \\
&  & 0  & 0   &   } $$ Consider the following push-out diagram:
$$\xymatrix{
& 0 \ar[d] &0 \ar[d] &   &    \\
0 \ar[r] & \Im f  \ar[d] \ar[r] &G_0
\ar[d]  \ar[r]  & M \ar@{=}[d] \ar[r]  &0  \\
0 \ar[r]  & B \ar[d]  \ar[r] &G
\ar[d] \ar[r] & M \ar[r]  &0  \\
& G_2\ar[d]\ar@{=}[r] &G_2 \ar[d] &   &    \\
&  0 & 0  &  & }$$ Because both $G_0$ and $G_2$ are Gorenstein
projective, $G$ is also Gorenstein projective by Lemma 2.1.
Connecting the middle rows in the above two diagrams, then we get
the first desired exact sequence. Since $G_0$ is Gorenstein
projective, there exists an exact sequence $0 \to G_3 \to Q \to G_0
\to 0$ in $\Mod R$ with $Q$ projective and $G_3$ Gorenstein
projective. Dually, taking pull-back, one gets the second desired
sequence. \hfill $\square$

\vspace{0.2cm}

For a positive integer $n$, recall that a module $A\in \Mod R$ is
called an {\it $n$-syzygy module} (of $M$) if there exists an exact
sequence $0\to A \to P_{n-1} \to \cdots \to P_1 \to P_0\to M \to 0$
in $\Mod R$ with all $P_i$ projective. Analogously, we give the
following

\vspace{0.2cm}

{\bf Definition 2.3.}  For a positive integer $n$, a module $A\in
\Mod R$ is a {\it Gorenstein $n$-syzygy module} (of $M$) if there
exists an exact sequence $0\to A \to G_{n-1} \to \cdots \to G_1 \to
G_0\to M \to 0$ in $\Mod R$ with all $G_i$ Gorenstein projective.

\vspace{0.2cm}

The following theorem is the main result in this section.

\vspace{0.2cm}

{\bf Theorem 2.4.} {\it Let $n$ be a positive integer and $0\to A\to
G_{n-1}\to G_{n-2}\to \cdots \to G_0\to M\to 0$ an exact sequence in
$\Mod R$ with all $G_i$ Gorenstein projective. Then we have the
following

(1) There exist exact sequences $0\to A\to P_{n-1}\to
P_{n-2}\to\cdots\to P_0\to N\to 0$ and $0\to M\to N\to G\to 0$ in
$\Mod R$ with all $P_i$ projective and $G$ Gorenstein projective. In
particular, a module in $\Mod R$ is an $n$-syzygy if and only if it
is a Gorenstein $n$-syzygy.

(2) There exist exact sequences $0\to B\to Q_{n-1}\to
Q_{n-2}\to\cdots\to Q_0\to M\to 0$ and $0\to H\to B\to A\to 0$ in
$\Mod R$ with all $Q_i$ projective and $H$ Gorenstein projective.}

\vspace{0.2cm}

{\it Proof.} (1) We proceed by induction on $n$. When $n=1$, it has
been proved in the proof of Proposition 2.2. Now suppose that $n\geq
2$ and we have an exact sequence: $$0\to A\to G_{n-1}\to G_{n-2}\to
\cdots \to G_0\to M\to 0$$ in $\Mod R$ with all $G_i$ Gorenstein
projective. Put $K=\Coker(G_{n-1}\to G_{n-2})$. By Proposition 2.2,
we get an exact sequence:
$$0\to A\to P_{n-1}\to G^{'}_{n-2}\to K\to 0$$
in $\Mod R$ with $P_{n-1}$ projective and $G^{'}_{n-2}$ Gorenstein
projective. Put $A^{'}=\Im(P_{n-1}\to G_{n-2}^{'})$. Then we get an
exact sequence:
$$0\to A^{'}\to G^{'}_{n-2}\to G_{n-3} \to \cdots \to G_0\to M\to 0$$
in $\Mod R$. So, by the induction hypothesis, we get the assertion.

(2) The proof is dual to that of (1), so we omit it. \hfill
$\square$

\vspace{0.2cm}

For a module $M\in \Mod R$, we use $\pd _{R}(M)$ to denote the
projective dimension of $M$.

\vspace{0.2cm}

{\bf Corollary 2.5.} ([CFH, Lemma 2.17]) {\it Let $M \in \Mod R$ and
$n$ a non-negative integer. If $\Gpd _{R}(M)=n$, then there exists
an exact sequence $0\to M \to N \to G \to 0$ in $\Mod R$ with $\pd
_{R}(N)=n$ and $G$ Gorenstein projective.}

\vspace{0.2cm}

{\it Proof.} Let $M \in \Mod R$ with $\Gpd _{R}(M)=n$. Then one uses
Theorem 2.4(1) with $A=0$ to get an exact sequence $0\to M\to N\to
G\to 0$ in $\Mod R$ with $\pd _{R}(N)\leq n$ and $G$ Gorenstein
projective. By Lemma 2.1, $\Gpd _{R}(N)=n$, and thus $\pd
_{R}(N)=n$. \hfill $\square$

\vspace{0.2cm}

By [H, Theorem 2.20], we have that $\Gpd _R(M) \leq n$ if and only
if there exists an exact sequence $0 \to G_{n} \to P_{n-1}\to \cdots
\to P_{1} \to P_{0} \to M \to 0$ in $\Mod R$ with all $P_{i}$
projective and $G_n$ Gorenstein projective. The following theorem
generalizes this result. In particular, the following theorem was
proved by Christensen and Iyengar in [CI, Theorem 3.1] when $R$ is a
commutative Noetherian ring.

\vspace{0.2cm}

{\bf Theorem 2.6.} {\it Let $M\in \Mod R$ and $n$ be a non-negative
integer. Then the following statements are equivalent.

(1) $\Gpd _R(M) \leq n$.

(2) For every non-negative integer $t$ such that $0 \leq t \leq n$,
there exists an exact sequence $0\to X_{n}\to \cdots \to X_1\to
X_0\to M\to 0$ in $\Mod R$ such that $X_t$ is Gorenstein projective
and $X_i$ is projective for $i \neq t$.}

{\it Proof.} $(2)\Rightarrow (1)$ It is trivial.

$(1)\Rightarrow (2)$ We proceed by induction on $n$. Suppose $\Gpd
_R(M) \leq 1$. Then there exists an exact sequence $0\to G_1\to
G_0\to M\to 0$ in $\Mod R$ with $G_0$ and $G_1$ Gorenstein
projective. By Proposition 2.2 with $A=0$, we get the exact
sequences $0\to P_1\to G'_0\to M\to 0$ and $0\to G'_1\to P_0\to M\to
0$ in $\Mod R$ with $P_0$, $P_1$ projective and $G'_0$ and $G'_1$
Gorenstein projective.

Now suppose $n\geq 2$. Then there exists an exact
sequence $0\to G_{n}\to \cdots \to G_1\to G_0\to M\to 0$ in $\Mod R$
with $G_i$ Gorenstein projective for any $1\leq i\leq n$. Set
$A=\Coker (G_3\to G_2)$. By applying Proposition 2.2 to the exact
sequence $0\to A\to G_1\to G_0\to M\to 0$, we get an exact sequence
$0\to G_{n}\to \cdots \to G_2\to G'_1\to P_0\to M\to 0$ in $\Mod R$
with $G'_1$ Gorenstein projective and $P_0$ projective. Set
$N=\Coker (G_2\to G'_1)$. Then we have $\Gpd _R(N) \leq n-1$. By the
induction hypothesis, there exists an exact sequence $0\to X_{n}\to
\cdots \to X_t\to \cdots \to X_1\to P_0\to M\to 0$ in $\Mod R$ such
that $P_0$ is projective and $X_t$ is Gorenstein projective and
$X_i$ is projective for $i\neq t$ and $1\leq t\leq n$.

Now we need only to prove (2) for $t=0$. Set $B=\Coker (G_2\to
G_1)$. By the induction hypothesis, we get an exact sequence $0\to
P_n\to \cdots\to P_3 \to P_2\to G'_1 \to B\to 0$ in $\Mod R$ with
$G'_1$ Gorenstein projective and $P_i$ projective for any $2\leq
i\leq n$. Set $C=\Coker (P_3\to P_2)$. Then by applying Proposition
2.2 to the exact sequence $0\to C\to G'_1\to G_0\to M\to 0$, we get
an exact sequence $0\to C\to P_1\to G'_0\to M\to 0$ in $\Mod R$ with
$P_0$ projective and $G'_0$ Gorenstein projective. Thus we obtain
the desired exact sequence $0\to P_n\to \cdots\to P_2 \to P_1\to
G'_0 \to M\to 0$. \hfill $\square$

\vspace{0.2cm}

Let $\mathscr{X}$ be a full subcategory of $\Mod R$. Recall from
[EJ2] that a morphism $f: X\to M$ in $\Mod R$ with $X \in
\mathscr{X}$ is called an {\it $\mathscr{X}$-precover} of $M$ if
$\Hom _{R}(X^{'}, X)
\xrightarrow{\Hom _{R}(X^{'},\ f)} \Hom _{R}(X^{'}, M) \to 0$ is exact for any
$X^{'}\in \mathscr{X}$. We use $\mathscr{GP}(R)$ to denote the full
subcategory of $\Mod R$ consisting of Gorenstein projective modules.
Let $M\in \Mod R$ with $\Gpd _R(M)=n<\infty$. Taking $t=0$ in
Theorem 2.6, one gets an exact sequence $0\to N\to G\to M\to 0$ in
$\Mod R$ with $G$ Gorenstein projective and $\pd _{R}(N)\leq n-1$.
It is easy to see that this exact sequence is a surjective
$\mathscr{GP}(R)$-precover of $M$ ([H, Theorem 2.10]).

\vspace{0.2cm}

{\bf Remark 2.7.} It is known that a module $A\in \Mod R$ is called
an {\it $n$-cosyzygy module} (of $M$) if there exists an exact
sequence $0\to M \to I^0 \to I^1 \to \cdots \to I^{n-1}\to A \to 0$
in $\Mod R$ with all $I^i$ injective. Recall from [EJ1] that a
module $E\in \Mod R$ is called {\it Gorenstein injective} if there
exists an exact sequence in $\Mod R$:
$$\cdots \to I_1 \to I_0 \to I^0 \to I^1\to \cdots ,$$
such that: (1) All $I_i$ and $I^i$ are injective; (2) After applying
the functor $\Hom _{R}(I, \ )$ the sequence is still exact for any
injective module $I \in \Mod R$; and (3) $E\cong \Im(I_0\to I^0)$.
We call $A$ a {\it Gorenstein $n$-cosyzygy module} (of $M$) if there
exists an exact sequence $0\to M \to E^0 \to E^1 \to \cdots \to
E^{n-1}\to A \to 0$ in $\Mod R$ with all $E^i$ Gorenstein injective.
We point out the dual versions on Gorenstein injectivity and
(Gorenstein) $n$-cosyzygy of all of the above results also hold true
by using a completely dual arguments.

\vspace{0.5cm}

\centerline{\bf 3. Gorenstein transpose}

\vspace{0.2cm}

In this section, $R$ is a left and right Noetherian ring and $\mod
R$ is the category of finitely generated left $R$-modules. For any
$A\in \mod R$, there exists a projective presentation in $\mod R$ :
$$P_{1} \buildrel {f}
\over \to P_{0} \to A \to 0.$$ Then we get an exact sequence $$0 \to
A^* \to P_{0}^* \buildrel {f^*} \over \longrightarrow P_{1}^* \to
\Coker f^* \to 0$$ in $\mod R^{op}$, where $(\ )^{*}=\Hom(\ , R)$.
Recall from [AB] that $\Coker f^*$ is called a {\it transpose} of
$A$, and denoted by $\Tr A$. We remark that the transpose of $A$
depends on the choice of the projective presentation of $A$, but it
is unique up to projective equivalence (see [AB]).

Analogously, we introduce the notion of Gorenstein transpose of
modules as follows. Let $A\in \mod R$. Then there exists a
Gorenstein projective presentation in $\mod R$: $$\pi: X_1 \buildrel
{g} \over \to X_0 \to A \to 0,$$ and we get an exact sequence:
$$0 \to A^* \to X_0^* \buildrel {g^*} \over
\longrightarrow X_1^* \to \Coker g^* \to 0$$ in $\mod R ^{op}$. We
call $\Coker g^*$ a {\it Gorenstein transpose} of $A$, and denote it
by $\Tr ^{\pi}_{G}A$. It is trivial that a transpose of $A$ is a
Gorenstein transpose of $A$, but the converse does not hold true in
general. For example, for a module $A$ in $\mod R$, if $A$ is Gorenstein projective
but not projective, then some Gorenstein transpose of $A$ is zero, and
any transpose of $A$ is Gorenstein projective (see Proposition 3.4(3) below) but
non-zero (otherwise, if a transpose of $A$ is zero, then $A$ is projective,
which is a contradiction).

Let $A\in \mod R$. Recall from [AB] that $A$ is said to have {\it Gorenstein
dimension zero} if $\Ext _{R}^i(A, R)=0=\Ext _{R^{op}}^i(\Tr A, R)$
for any $i \geq 1$. It is easy to see that if $A$ has Gorenstein
dimension zero, then so does $A^{*}$. In addition, it is well known
that $A$ has Gorenstein dimension zero if and only if it is
Gorenstein projective. Let $\sigma _{A}: A \to A^{**}$ defined via
$\sigma _{A}(x)(f)=f(x)$ for any $x\in A$ and $f\in A^{*}$ be the
canonical evaluation homomorphism.
Recall that a module $A\in \mod R$ is called {\it torsionless}
(resp. {\it reflexive}) if $\sigma _A$ is a monomorphism (resp. an
isomorphism)

The following result establishes a relation between a Gorenstein
transpose of a module with a transpose of the same module.

\vspace {0.2cm}

{\bf Theorem 3.1.} {\it Let $M\in \mod R^{op}$ and $A\in \mod R$.
Then $M$ is a Gorenstein transpose of $A$ if and only if $M$ can be
embedded into a transpose $\Tr A$ of $A$ with the cokernel
Gorenstein projective, that is, there exists an exact sequence $0\to
M\to \Tr A\to H\to 0$ in $\mod R^{op}$ with $H$ Gorenstein
projective.}

\vspace {0.2cm}

{\it Proof.} We first prove the necessity. Assume that $M(\cong \Tr
^{\pi}_G A)$ is a Gorentein transpose of $A$. Then there exists an exact sequence $\pi: X_1 \buildrel {g}
\over \to X_0 \to A \to 0$ in $\mod R$ with $X_0$ and $X_1$
Gorenstein projective such that $\Tr ^{\pi}_G A=\Coker g^{*}$.
So there exists an exact sequence $0 \to H^{'}_1 \to P_0^{'} \to X_0 \to 0$
in $\mod R$ with $P_0^{'}$ projective and $H^{'}_1$
Gorenstein projective. Let $K_1=\Im g$ and $g=i\alpha$
be the natural epic-monic decomposition of $g$. Then we have the following
pull-back diagram:

$$\xymatrix{
&  & 0 \ar[d]  & 0 \ar[d] &  & \\
&  &  H^{'}_1 \ar[d] \ar@{=}[r] & H^{'}_1 \ar[d]  &  \\
& 0 \ar[r] & K^{'}_1\ar[d] \ar[r] & P_0^{'} \ar[d]  \ar[r]  & A \ar@{=}[d] \ar[r] & 0 \\
& 0 \ar[r] & K_1\ar[d]\ar[r]^{i} & X_0 \ar[d]  \ar[r]  & A \ar[r] & 0 \\
&  & 0   & 0   &   & }$$ Now consider the following
pull-back diagram:
$$\xymatrix{
&  &  & 0 \ar[d]  & 0 \ar[d] &   \\
&  &  &H^{'}_1 \ar[d] \ar@{=}[r] & H^{'}_1 \ar[d] & \\
& 0 \ar[r] & K_2\ar@{=}[d] \ar[r]  &G \ar[d]  \ar[r]  & K^{'}_1 \ar[d] \ar[r] & 0 \\
& 0 \ar[r] & K_2\ar[r] & X_1 \ar[d]  \ar[r]^{\alpha}  & K_1 \ar[r]\ar[d] & 0 \\
&  &   & 0   &  0 & }$$ where $K_2=\Ker g$. Because both $X_1$ and $H^{'}_1$
are Gorenstein projective, $G$ is Gorenstein projective by Lemma 2.1.
So there exists an exact sequence $0 \to G_1 \to P_0 \to G \to 0$
in $\mod R$ with $P_0$ projective and $G_1$
Gorenstein projective. Consider the following pull-back diagram:

$$\xymatrix{
&  & 0 \ar[d]  & 0 \ar[d] &  & \\
&  &  G_1 \ar[d] \ar@{=}[r] & G_1 \ar[d]  &  \\
& 0 \ar[r] & K^{'}_2\ar[d]^{\beta} \ar[r]  &P_0 \ar[d]  \ar[r]  & K^{'}_1 \ar@{=}[d] \ar[r] & 0 \\
& 0 \ar[r] & K_2\ar[d]\ar[r] & G \ar[d]  \ar[r]  & K^{'}_1 \ar[r] & 0 \\
&  & 0   & 0   &   & }$$
So we get the following commutative diagram with exact rows:
$$\xymatrix{
& 0 \ar[r] & K^{'}_2\ar[d]^{\beta} \ar[r]  &P_0 \ar[d]  \ar[r] & K^{'}_1 \ar@{=}[d] \ar[r] & 0 \\
& 0 \ar[r] & K_2\ar@{=}[d]\ar[r] & G \ar[d]  \ar[r]  & K^{'}_1 \ar[r]\ar[d] & 0 \\
& 0 \ar[r] & K_2\ar[r] & X_1  \ar[r]^{\alpha}  & K_1 \ar[r] & 0}
$$
It yields the following commutative diagram with exact columns and rows:
$$\xymatrix{
&  & 0\ar[d] & 0 \ar[d]  & 0 \ar[d] &   \\
&  & \Ker\beta\ar[d] &H_1 \ar[d]  & H^{'}_1 \ar[d] &   \\
& 0 \ar[r] & K^{'}_2\ar[d]^{\beta} \ar[r]  &P_0 \ar[d]  \ar[r]  & K^{'}_1 \ar[d] \ar[r] & 0 \\
& 0 \ar[r] & K_2\ar[r]\ar[d] & X_1 \ar[d]  \ar[r]^{\alpha}  & K_1 \ar[r]\ar[d] & 0 \\
&  &  0 & 0   &  0 & }$$ where $H_1=\Ker(P_0\to X_1)$.
By the snake lemma, we get the exact
sequence $0\to \Ker \beta\to H_1\buildrel {h} \over\to H^{'}_1\to
0$. By Lemma 2.1, $H_1$ is
Gorenstein projective and hence $\Ker \beta$ is also Gorenstein
projective. Combining the above diagram with the first one in this proof,
we get the following commutative diagram with
exact columns and rows:

$$\xymatrix{
& &  0\ar[d]& 0 \ar[d]  & 0 \ar[d] &  & \\
&0\ar[r]&\Ker \beta\ar[d]\ar[r]  & H_1 \ar[d] \ar[r]^{h} & H^{'}_1 \ar[d]\ar[r] & 0 & \\
&0\ar[r] &K^{'}_2 \ar[d]^{\beta}\ar[r]& P_0\ar[d] \ar[r] &P_0^{'} \ar[d]  \ar[r]  & A \ar@{=}[d] \ar[r] & 0 \\
&0\ar[r] &K_2\ar[r]\ar[d] & X_1\ar[d]\ar[r]^{g} & X_0 \ar[d]  \ar[r]  & A \ar[r] & 0 \\
&  & 0   & 0   &  0 & }$$ By applying the
functor $(\ )^*$ to the above diagram, we get the following
commutative diagram with exact columns and rows:
$$\xymatrix{
&  & 0   & 0  &  & \\
& & H_1^* \ar[u] & {H^{'}_1}^* \ar[u] \ar[l]_{h^*}  & \ar[l] 0 & \\
& & P_0^* \ar[u]  & {P_0^{'}}^* \ar[u]  \ar[l]  & A^*  \ar[l] & \ar[l] 0 \\
& & X_1^* \ar[u] & X_0^* \ar[u] \ar[l]_{g^*}  & A^* \ar@{=}[u]\ar[l] & \ar[l] 0 \\
&  & 0 \ar[u] & 0 \ar[u]  &   & }$$ By the snake lemma, we get an
exact sequence: $$0\to \Tr ^{\pi}_{G}A(=\Coker g^*)\to \Tr A\to
\Coker h^*\to 0$$ in $\mod R^{op}$ with $\Coker h^*(\cong
(\Ker h)^*\cong(\Ker \beta)^*)$ Gorenstein projective.

We next prove the sufficiency. Let $P_1\buildrel {f} \over \to P_0\to A\to 0$ be a
projective presentation of $A$ in $\mod R$. Then we have the
following pull-back diagram:
$$\xymatrix{
&  & 0 \ar[d]  & 0 \ar[d] &  &  \\
&  & A^* \ar[d] \ar@{=}[r] & A^* \ar[d] &  & \\
&  & P_0^* \ar[d]^{h} \ar@{=}[r] & P_0^* \ar[d]^{f^*} &  & \\
&0\ar[r] & K\ar[d]^{g} \ar[r] & P_1^* \ar[d] \ar[r]  & H \ar@{=}[d] \ar[r] & 0 \\
& 0\ar[r]& M\ar[d]\ar[r] & \Tr A \ar[d]  \ar[r]  & H \ar[r] & 0 \\
&  & 0   & 0   &   & .}
$$ Because $H$ is Gorenstein projective and $P_1^*$ is projective, $K$
is Gorenstein projective by Lemma 2.1. Again because $H$ is
Gorenstein projective, by applying the functor $(\ )^*$ to the above
commutative diagram, we get the following commutative diagram with
exact columns and rows:

$$\xymatrix{
&  &  & 0 \ar[d] &  0 \ar[d]&  &\\
& 0\ar[r] & H^*\ar[r]\ar@{=}[d] & (\Tr A)^*\ar[r] \ar[d] &  M^* \ar[r]\ar[d]^{g^*}& 0 &\\
& 0\ar[r] & H^*\ar[r]& P_1^{**}\ar[r] \ar[d]^{f^{**}} &  K^* \ar[r]\ar[d]^{h^*}& 0 &\\
&  &  & P_0^{**} \ar@{=}[r]\ar[d] &  P_0^{**} &  & \\
&  &  & A \ar[d] &  &  &\\
&  &  & 0 &  &  &}
$$
By the snake Lemma, we have $\Im h^*\cong \Im f^{**}$. Thus we get
$\Coker h^*=P_0^{**}/\Im h^*\cong P_0^{**}/\Im f^{**} \cong A$, and
therefore we get a Gorenstein projective presentation of $A$ in
$\mod R$:
$$K^*\buildrel {h^*} \over \longrightarrow P_0^{**}\to A\to 0.$$
Because both $K$ and
$P_0^*$ are reflexive, we get an exact sequence $0\to A^*\to
P_0^{***}\buildrel {h^{**}} \over \longrightarrow K^{**}\to M\to 0 $
in $\mod R^{op}$ and $M$ is a Gorenstein transpose of $A$. \hfill
$\square$

\vspace{0.2cm}

As a consequence of Theorem 3.1, we get the following

\vspace{0.2cm}

{\bf Corollary 3.2.} {\it Let $A\in \mod R$. Then for any Gorenstein
projective module $H\in \mod R^{op}$ and any transpose $\Tr A$ of
$A$, $H \oplus \Tr A$ is a Gorenstein transpose of $A$.}

\vspace {0.2cm}

{\it Proof.} Assume that $H\in \mod R^{op}$ is a Gorenstein
projective module. Then there exists an exact sequence $0\to H \to P
\to H' \to 0$ in $\mod R^{op}$ with $P$ projective and $H'$
Gorenstein projective, which induces an exact sequence $0\to H\oplus
\Tr A \to P\oplus \Tr A \to H' \to 0$. Because $P\oplus \Tr A$ is
again a transpose of $A$, $H\oplus \Tr A$ is a Gorenstein transpose
of $A$ by Theorem 3.1.  \hfill $\square$

\vspace{0.2cm}

It is clear that the Gorenstein transpose of a module $A$ in $\mod R$
depends on the choice of the Gorenstein projective presentation of
$A$. Corollary 3.2 provides a method to construct a Gorenstein transpose
of a module from a transpose of the same module. It is interesting to ask the following

\vspace{0.2cm}

{\bf Question 3.3.} Is any Gorenstein transpose obtained in this way?

\vspace{0.2cm}

If the answer to this question is positive, then we can conclude that
the Gorenstein transpose of a module is unique up to Gorenstein
projective equivalence.

Let $A\in \mod R$. By [A, Proposition 6.3] (or [AB,
Proposition 2.6]), there exists an exact sequence:
$$0 \to \Ext _{R^{op}}^{1}(\Tr A, R) \to A \buildrel {\sigma _{A}}
\over \longrightarrow A^{**} \to \Ext _{R^{op}}^{2}(\Tr A, R) \to 0
\eqno{(*)}$$ in $\mod R$. For a positive integer $n$, recall from [AB] that $A$
is called {\it $n$-torsionfree} if $\Ext _{R^{op}}^i(\Tr A,
R)=0$ for any $1\leq i\leq n$. From the exact sequence $(*)$, it is
easy to see that $A$ is torsionless (resp.
reflexive) if and only if it is 1-torsionfree (resp. 2-torsionfree).

The following result shows that some homological properties of any
Gorenstein transpose and any transpose of a given module are identical.

\vspace{0.2cm}

{\bf Proposition 3.4.} {\it Let $A\in \mod R$. Then for any
Gorenstein transpose $\Tr_G^{\pi} A$ and
any transpose $\Tr A$ of $A$, we have

(1) $\Ext _{R^{op}}^i(\Tr_G^{\pi} A, R)\cong \Ext _{R^{op}}^{i}(\Tr A, R)$
for any $i \geq 1$.

(2) For any $n\geq 1$, $\Tr_G^{\pi} A$ is $n$-torsionfree if and only if so is $\Tr A$.

(3) Some Gorenstein transpose of $A$ is zero if and only if $A$ is
Gorenstein projective, if and only if any (Gorenstein) transpose of
$A$ is Gorenstein projective.

(4) $\Gpd _{R^{op}}(\Tr_G^{\pi} A)= \Gpd _{R^{op}}(\Tr A)$.}

\vspace{0.2cm}

{\it Proof.} (1) It is an immediate consequence of Theorem 3.1.

(2) Let $\Tr_G^{\pi} A$ be any Gorenstein transpose of $A$. By
Theorem 3.1, there exists a transpose $\Tr A$ of $A$ satisfying the
exact sequence $0\to \Tr_G^{\pi} A\to \Tr A\to H\to 0$ in $\mod
R^{op}$ with $H$ Gorenstein projective.

If $\Ext_R^1(\Tr (\Tr A),R)=0$, then $\Tr A$ is torsionless. So
$\Tr_G^{\pi} A$ is also torsionless and $\Ext_R^1(\Tr (\Tr_G^\pi A),
R)=0$. Because $H$ is Gorenstein projective, we get an exact
sequence $0\to \Tr H \to \Tr (\Tr A)\to \Tr (\Tr_G^\pi A)\to 0$ in
$\mod R$ with $\Tr H$ Gorenstein projective. So we have that
$\Ext_R^i(\Tr (\Tr_G^\pi A), R)\cong \Ext_R^i(\Tr (\Tr A),R)$ for
any $i\geq 2$, and $\Ext_R^1(\Tr (\Tr_G^\pi A), R)\to \Ext_R^1(\Tr
(\Tr A),R) \to 0$ is exact. So for any $i\geq 1$, $\Ext_R^i(\Tr
(\Tr_G^\pi A), R)=0$ if and only if $\Ext_R^i(\Tr (\Tr A),R)=0$, and
thus we conclude that for any $n \geq 1$, $\Tr_G^{\pi} A$ is
$n$-torsionfree if and only if so is $\Tr A$.

(3) Because $A$ is a (Gorenstein) transpose of any (Gorenstein)
transpose of $A$, it is not difficult to verify the assertion by (1)
and (2).

(4) Let $\Tr_G^{\pi} A$ be any Gorenstein transpose of $A$. If
$\Tr_G^{\pi} A=0$, then the assertion follows from (3). Now suppose
$\Tr_G^{\pi} A\neq 0$. By Theorem 3.1, there exists a transpose $\Tr
A$ of $A$ satisfying the exact sequence $0\to \Tr_G^{\pi} A\to \Tr
A\to H\to 0$ in $\mod R^{op}$ with $H$ Gorenstein projective. Then
we have that $\Gpd _{R^{op}}(\Tr_G^{\pi} A)= \Gpd _{R^{op}}(\Tr A)$
by Lemma 2.1. \hfill{$\square$}

\vspace{0.2cm}

Let $A \in \mod R$. By Proposition 3.4(1), we have that $A$ is
$n$-torsionfree if and only if $\Ext _{R^{op}}^i(\Tr_G^{\pi} A,
R)=0$ for any (or some) Gorenstein transpose $\Tr_G^{\pi} A$ of $A$
and $1\leq i\leq n$. On the other hand, also by Proposition 3.4(1),
we get a Gorenstein version of the
formula $(*)$ as follows. For any Gorenstein transpose $\Tr ^{\pi}_GA$ of
$A$, we have the following exact sequence:
$$0 \to \Ext _{R^{op}}^{1}(\Tr ^{\pi}_{G}A, R) \to A \buildrel {\sigma _{A}}
\over \longrightarrow A^{**} \to \Ext _{R^{op}}^{2}(\Tr ^{\pi}_{G}A,
R) \to 0$$ in $\mod R$. It is easy to see that $A$ is a Gorenstein transpose of
$\Tr ^{\pi}_{G}A$. So we also get the following exact sequence:
$$0 \to \Ext_R^{1}(A, R) \to
\Tr_G^{\pi} A \buildrel {\sigma _{\Tr ^{\pi}_{G}A}} \over \longrightarrow
(\Tr_G^{\pi} A)^{**} \to \Ext_R^{2}(A, R) \to 0$$ in $\mod R^{op}$.

The following result shows that any double Gorenstein transpose of
$A$ shares some homological properties of $A$.

\vspace{0.2cm}

{\bf Corollary 3.5.} {\it Let $A\in \mod R$. Then for any Gorenstein
transpose $\Tr_G^{\pi} A$ of $A$ and any Gorenstein transpose
$\Tr_G^{\pi '}(\Tr_G^{\pi} A)$ of $\Tr_G^{\pi} A$, we have

(1) $\Ext_R^i(\Tr_G^{\pi '}(\Tr_G^{\pi} A), R)\cong \Ext _R^{i}(A,
R)$ for any $i \geq 1$.

(2) For any $n\geq 1$, $\Tr_G^{\pi '}(\Tr_G^{\pi} A)$ is
$n$-torsionfree if and only if so is $A$.

(3) $\Gpd _R(\Tr_G^{\pi '}(\Tr_G^{\pi} A))=\Gpd _R(A)$.}

\vspace{0.2cm}

{\it Proof.} Note that $A$ is a Gorenstein transpose of any
Gorenstein transpose $\Tr_G^{\pi} A$ of $A$. So all of the
assertions follow from Proposition 3.4.  \hfill{$\square$}

\vspace{0.2cm}

Note that a transpose of a module is a special Gorenstein transpose
of the same module. The following result shows that a module with
Gorenstein projective dimension $n$ is a double Gorenstein transpose
of a module with projective dimension $n$.

\vspace {0.2cm}

{\bf Proposition 3.6.} {\it Let $A\in \mod R$ and $n$ be a
non-negative integer. Then $\Gpd _R(A)=n$ if and only if there
exists a module $B\in \mod R$ with $\pd _{R}(B)=n$ such that $A$ is
a Gorenstein transpose of some transpose $\Tr B$ of $B$ (that is,
$A=\Tr ^{\pi}_G(\Tr B)$, where $\Tr ^{\pi}_G(\Tr B)$ is a Gorenstein
transpose of some transpose $\Tr B$ of $B$).}

\vspace{0.2cm}

{\it Proof.} Assume that $A\in \mod R$ with $\Gpd _R(A)=n$. By
Corollary 2.5, there exists an exact sequence $0\to A \to B \to H\to
0$ in $\mod R$ with $\pd _{R}(B)=n$ and $H$ Gorenstein projective.
Note that $B$ is a transpose of some transpose $\Tr B$ of $B$. By
Theorem 3.1, $A$ is a Gorenstein transpose of $\Tr B$.

Conversely, if $A$ is a Gorenstein transpose of some transpose $\Tr
B$ of a module $B\in \mod R$ with $\pd _{R}(B)=n$, then $\Gpd
_R(A)=\Gpd _R(B)=\pd _{R}(B)=n$ by Corollary 3.5. \hfill $\square$

\vspace{0.5cm}

{\bf Acknowledgements.} This research was partially supported by the
Specialized Research Fund for the Doctoral Program of Higher
Education, NSFC (Grant No. 10771095) and NSF
of Jiangsu Province of China (Grant Nos. BK2010047, BK2010007). The authors
thank Lars Winther Christensen and the referee for the useful suggestions.


\vspace{0.5cm}

\begin{thebibliography}{101}

\bibitem[A]{A1} M. Auslander, {\it Coherent functors}, in: Proc. of the Conf. on Categorial Algebra, La Jolla, 
1965, Springer-Verlag, Berlin, 1966, pp.189--231.

\bibitem[AB]{A2} M. Auslander and M. Bridger, {\it Stabe Module Theory},
Memoirs Amer. Math. Soc. {\bf 94}, Amer. Math. Soc., Providence, RI,
1969.

\bibitem[C]{A3} L.W. Christensen, Gorenstein Dimension, Lect. Notes in Math. {\bf 1747},
Springer-Verlag, Berlin, 2000.

\bibitem[CFH]{A4} L.W. Christensen, A. Frankild and H. Holm, {\it On Gorenstein projective,
injective and flat dimensions--A functorial description with
applications}, J. Algebra, {\bf 302} (2006), 231--279.

\bibitem[CI]{A5} L.W. Christensen and S. Iyengar, {\it Gorenstein
dimension of modules over homomorphisms}, J. Pure Appl. Algebra,
{\bf 208} (2007), 177--188.

\bibitem[EJ1]{A6} E.E. Enochs and O.M.G. Jenda, {\it Gorenstein injective and
projective modules}, Math. Z., {\bf 220} (1995), 611--633.

\bibitem[EJ2]{A7} E.E. Enochs and O.M.G. Jenda, Relative
Homological Algebra, De Gruyter Exp. in Math. {\bf 30}, Walter de
Gruyter, Berlin, New York, 2000.

\bibitem[H]{A8} H. Holm, {\it Gorenstein homological dimensions}, J.
Pure Appl. Algebra, {\bf 189} (2004), 167--193.

\end{thebibliography}
\end{document}